\documentclass{article}

\newcommand{\be}{\begin{equation}}
\newcommand{\ee}{\end{equation}}
\newcommand{\bea}{\begin{eqnarray}}
\newcommand{\eea}{\end{eqnarray}}
\newcommand{\barray}{\begin{array}}
\newcommand{\earray}{\end{array}}
\newcommand{\pa}{\partial}
\newcommand{\nn}{\nonumber}
\newcommand{\bitem}{\begin{itemize}}
\newcommand{\eitem}{\end{itemize}}
\newtheorem{teo}{Theorem}[section]
\newcommand{\bt}{\begin{teo}}
\newcommand{\et}{\end{teo}}
\newtheorem{Def}{Definition}[section]
\newcommand{\bd}{\begin{Def}}
\newcommand{\ed}{\end{Def}}
\newtheorem{lem}{Lemma}[section]
\newcommand{\bl}{\begin{lem}}
\newcommand{\el}{\end{lem}}
\newtheorem{prop}{Proposition}[section]
\newcommand{\bp}{\begin{prop}}
\newcommand{\ep}{\end{prop}}
\newtheorem{cor}{Corollary}[section]
\newcommand{\bc}{\begin{cor}}
\newcommand{\ec}{\end{cor}}
\newtheorem{ex}{Example}[section]
\newcommand{\bex}{\begin{ex}}
\newcommand{\eex}{\end{ex}}
\newtheorem{rem}{Remark}[section]
\newcommand{\br}{\begin{rem}}
\newcommand{\er}{\end{rem}}

\catcode `\@=11
\@addtoreset{equation}{section}

\begin{document}

\begin{center}
{\Large \textbf{Liouville canonical form for
compatible \\ nonlocal Poisson brackets of
hydrodynamic \\ type, and integrable
hierarchies\footnote{This work was supported by
the Alexander von Humboldt Foundation (Germany),
the Russian Foundation for Basic Research
(grant No. 99--01--00010) and the INTAS
(grant No. 99--1782).}}}
\end{center}

\bigskip
\bigskip

\centerline{\large {O. I. Mokhov}}

\bigskip
\medskip

\section{Introduction}

In the present paper, we solve the problem of reducing to
the simplest and convenient for our purposes,
``canonical'' form for an arbitrary pair of
compatible nonlocal Poisson brackets of hydrodynamic type
generated by metrics of constant Riemannian curvature
(compatible Mokhov--Ferapontov brackets \cite{1})
in order to get an effective construction of the integrable hierarchies
related to all these compatible Poisson brackets.
As was shown in \cite{2}, \cite{3}
(see also \cite{4}), compatible Mokhov--Ferapontov
brackets are described by a consistent nonlinear system of
equations integrable by the method of
inverse scattering problem (the case of flat metrics
see in \cite{5}--\cite{7}). But the problem
of an effective construction of the corresponding integrable
hierarchies in this case, what is the main purpose of the
present paper, requires a different approach to
the description of these compatible brackets.
In this paper, for an arbitrary solution of the
integrable system of equations describing the compatible
brackets under consideration, that is, for an arbitrary
pair of these compatible brackets, integrable
bi-Hamiltonian systems of hydrodynamic type possessing
this pair of compatible nonlocal Poisson brackets of
hydrodynamic type are constructed in an explicit form.
For the case of the Dubrovin--Novikov brackets \cite{8}
(the local Poisson brackets of hydrodynamic type),
this problem was considered and completely solved in
the present author's works \cite{9}, \cite{10}.

In \cite{1}
the nonlocal Poisson brackets of hydrodynamic type
which have the following form (the Mokhov--Ferapontov brackets):
\be
\{ I,J \} = \int {\delta I \over \delta u^i(x)}
\left ( g^{ij}(u(x)) {d \over dx} + b^{ij}_k (u(x))\, u^k_x
+ K u^i_x \left ( {d \over dx} \right )^{-1} u^j_x \right )
{\delta J \over \delta u^j(x)} dx,
\label{nonl}
\ee
where $I [u]$ and $J [u]$ are arbitrary
functionals on the space of functions (fields)
$u^i (x),$ $1 \leq i \leq N,$
of single independent variable $x$, $u = (u^1,..., u^N)$ are
local coordinates on a certain given smooth $N$-dimensional
manifold $M$, the coefficients $g^{ij} (u)$ and $b^{ij}_k (u)$
of the bracket (\ref{nonl}) are smooth functions of
local coordinates,
$K$ is an arbitrary constant,
were introduced and studied.

The form of bracket (\ref{nonl}) is invariant
with respect to local changes of coordinates.
A bracket of the form (\ref{nonl})
is called {\it nondegenerate} if $\det (g^{ij} (u)) \not\equiv 0$.
If $\det (g^{ij} (u)) \not\equiv 0$, then bracket (\ref{nonl})
is a Poisson bracket if and only if
$g^{ij} (u)$ is an arbitrary pseudo-Riemannian
contravariant metric of constant Riemannian curvature $K$,
$b^{ij}_k (u) = - g^{is} (u) \Gamma ^j_{sk} (u),$ where
$\Gamma^j_{sk} (u)$ is the Riemannian connection
generated by the metric $g^{ij} (u)$
(the Levi--Civita connection) \cite{1}
(note that
the coefficients
$g^{ij} (u)$ and $b^{ij}_k (u)$ of bracket
(\ref{nonl}) are transformed as
corresponding differential-geometric objects
under local changes of coordinates: a contravariant
metric
$g^{ij} (u)$ and a contravariant connection $b^{ij}_k (u)
= - g^{is} (u) \Gamma^j_{sk} (u)$
respectively, $K$ is an invariant).
For $K= 0$  we have the local Poisson brackets
of hydrodynamic type
(the Dubrovin--Novikov brackets \cite{8}).

Recall that Poisson brackets
are called {\it compatible} if their arbitrary
linear combination is also a Poisson bracket
(Magri, \cite{11}).

\section{Compatible nonlocal Poisson brackets \\
of hydrodynamic type}

\bl \label{locnon}
In the classification problem
for an arbitrary pair of compatible
nonlocal Poisson brackets of the form (\ref{nonl}),
one can always consider one of the Poisson brackets as
local without loss of generality.
\el

Actually, if two compatible nonlocal Poisson brackets
of the form (\ref{nonl})
$\{ I, J \}_0$ (with a corresponding constant $K_0$ in the nonlocal
term) and $\{ I, J \}_1$ (with a constant $K_1$)
are linear independent, then in the pencil of these
Poisson brackets, that is, among the Poisson brackets
 $\lambda_0 \{ I, J \}_0 +
\lambda_1 \{ I, J \}_1$, where $\lambda_0$ and $\lambda_1$
are arbitrary constants, there is necessarily
a nonzero local Poisson bracket
(for this bracket $\lambda_0 K_0 + \lambda_1 K_1 = 0$),
which can be taken as one of the generators for all the considered
pencil of compatible Poisson brackets.

Consider the problem of compatibility
for a pair of
nonlocal and local Poisson brackets of hydrodynamic type
\be
\{ I,J \}_1 = \int {\delta I \over \delta u^i(x)}
\left ( g^{ij}_1 (u(x)) {d \over dx} + b^{ij}_{1, k} (u(x))\, u^k_x
+ K_1 u^i_x \left ( {d \over dx} \right )^{-1} u^j_x \right )
{\delta J \over \delta u^j(x)} dx
\label{nonloc}
\ee
and
\be
\{ I,J \}_2 = \int {\delta I \over \delta u^i(x)}
\left ( g^{ij}_2 (u(x)) {d \over dx} + b^{ij}_{2, k} (u(x))\, u^k_x
 \right )
{\delta J \over \delta u^j(x)} dx,
\label{lok}
\ee
that is,
the condition that for any constant $\lambda$
the bracket
\be
\{ I,J \} = \{ I, J \}_1 + \lambda \{ I, J \}_2
\label{comp}
\ee
is a Poisson bracket ({\it a pencil of Poisson brackets}).

We assume further that the local bracket $\{ I, J \}_2$
is nondegenerate, that is,
$\det g^{ij}_2 (u) \not\equiv 0$,
but we do not impose any additional conditions
on the bracket $\{ I, J \}_1$, that is, generally speaking,
it may be degenerate.
The bracket (\ref{comp}) can be degenerate, therefore
we need general relations for
coefficients of an arbitrary bracket of the form (\ref{nonl})
which are equivalent to the condition that
the bracket (\ref{nonl}) is a Poisson bracket.
These general relations (without the assumption of
nondegeneracy) were obtained in the present author's work
\cite{12}
(see also \cite{13}):
\be
g^{ij} (u) = g^{ji} (u), \label{s1}
\ee
\be
{\pa g^{ij} \over \pa u^k} = b^{ij}_k (u) + b^{ji}_k (u),\label{s2}
\ee
\be
g^{is} (u) b^{jr}_s (u) = g^{js} (u) b^{ir}_s (u), \label{s3}
\ee
\be
g^{is} (u) \left ( {\pa b^{jr}_s \over \pa u^k}  -
{\pa b^{jr}_k \over \pa u^s} \right )
+ b^{ij}_s (u) b^{sr}_k (u) - b^{ir}_s (u) b^{sj}_k (u) =
K ( g^{ir} (u) \delta^j_k -
g^{ij} (u) \delta^r_k ),\label{s4}
\ee
\bea
&&
\sum_{(i, j, r)} \left [
b^{si}_p (u) \left (
{\pa b^{jr}_k \over \pa u^s} -
{\pa b^{jr}_s \over \pa u^k} \right ) + K (b^{ij}_k (u) - b^{ji}_k (u))
\delta^r_p +\right. \nn\\
&&
\left. b^{si}_k (u) \left (
{\pa b^{jr}_p \over \pa u^s} -
{\pa b^{jr}_s \over \pa u^p} \right ) + K (b^{ij}_p (u) - b^{ji}_p (u))
\delta^r_k \right ]  = 0, \label{s5}
\eea
where $\sum_{(i, j, r)}$ means summation over all cyclic
permutations of the indices $i, j, r$.

\section{Canonical form for compatible
pairs of brackets}

According to the Dubrovin--Novikov theorem
\cite{8},
for any nondegenerate local Poisson bracket
of hydrodynamic type $\{ I, J \}_2$,
there always exist local coordinates $u^1,...,u^N$
(flat coordinates of the metric $g^{ij}_2 (u)$)
in which this bracket is constant, that is,
$g^{ij}_2 (u) = \eta^{ij} = {\rm\ const}$, $b^{ij}_{2, k} (u) =
\Gamma^i_{2, jk} (u) = 0$.
Thus we can choose flat coordinates of the metric
$g^{ij}_2 (u)$
and further on consider that the Poisson bracket
$\{ I, J \}_2$ is constant and has the form
\be
\{ I,J \}_2 = \int {\delta I \over \delta u^i(x)}
 \eta^{ij} {d \over dx}
{\delta J \over \delta u^j(x)} dx,
\label{lokconst}
\ee
where $\eta^{ij} = \eta^{ji}$, $\eta^{ij} = {\rm\ const},$
$\det \eta^{ij} \neq 0,$ $\eta^{is} \eta_{sj} = \delta^i_j.$
In the sequel, in the considered flat coordinates,
we shall also
use the covariant metric $\eta_{ij}$ which is
inverse to the contravariant metric
$\eta^{ij}$.

\bt \label{t1}
An arbitrary nonlocal Poisson bracket $\{ I, J \}_1$
of the form (\ref{nonloc})
(may be degenerate) is compatible with the
constant Poisson bracket
(\ref{lokconst}) if and only if
it has the form
\bea
&&
\ \{ I,J \}_1 = \int {\delta I \over \delta u^i(x)}
\left ( \left [ \eta^{is} {\pa H^j \over
\pa u^s} + \eta^{js} {\pa H^i \over \pa u^s}
- K_1 u^i u^j \right ] {d \over dx} + \right.\nn\\
&&
\left. \left [
\eta^{is} {\pa^2 H^j \over \pa u^s \pa u^k} - K_1 \delta^i_k u^j
\right ]\, u^k_x
+ K_1 u^i_x \left ( {d \over dx} \right )^{-1} u^j_x \right )
{\delta J \over \delta u^j(x)} dx,
\label{nonloc3}
\eea
where $H^i (u),$ $1 \leq i \leq N,$ are smooth
functions defined in a certain domain of local coordinates.
\et

In the flat case of compatible
Dubrovin--Novikov brackets
($K_1 = 0$), the corresponding, necessary for our
purposes, statement was formulated by the
present author in  \cite{14}, \cite{16}
(see also the conditions on flat pencils of metrics in \cite{18}).

{\it Proof.}
It follows from relations (\ref{s1})--(\ref{s5})
that, in the considered local coordinates,
the conditions of compatibility for the Poisson brackets
$\{ I, J \}_1$
and $\{ I, J \}_2$ or, in other words, the conditions that
the bracket
\bea
&&
\ \{ I,J \} = \int {\delta I \over \delta u^i(x)}
\left ( ( g^{ij}_1 (u(x)) +
\lambda \eta^{ij} ){d \over dx} + \right.\nn\\
&&
\left. b^{ij}_{1, k} (u(x))\, u^k_x
+ K_1 u^i_x \left ( {d \over dx} \right )^{-1} u^j_x \right )
{\delta J \over \delta u^j(x)} dx
\label{nonloc2}
\eea
is a Poisson bracket for all values of the parameter $\lambda$
have the form
\be
\eta^{is}b^{jr}_{1, s} (u) = \eta^{js} b^{ir}_{1, s} (u), \label{c1}
\ee
\be
{\pa b^{jr}_{1, s} \over \pa u^k}  -
{\pa b^{jr}_{1, k} \over \pa u^s}  =
K_1 ( \delta^r_s \delta^j_k  -
\delta^j_s \delta^r_k ).\label{c2}
\ee

Let us define the functions $A^{ij}_k (u)$
by the relations
\be
A^{ij}_k (u) = b^{ij}_{1, k} (u) - K_1 \delta^j_k u^i. \label{for1}
\ee
It follows from formula (\ref{c2}) that
$${\pa A^{jr}_s \over \pa u^k}  -
{\pa A^{jr}_k \over \pa u^s}  = 0,$$
that is, by the Poincar\'e lemma, there locally exist functions
$P^{ij} (u)$ such that
$$A^{ij}_k (u) = {\pa P^{ij} \over \pa u^k},$$
and we derive from (\ref{for1}) a necessary
expression for the coefficient $b^{ij}_{1, k} (u)$:
$$b^{ij}_{1, k} (u) =
{\pa P^{ij} \over \pa u^k} + K_1 \delta^j_k u^i.$$
Let us find the corresponding
expression for the metric $g^{ij}_1 (u)$.
From relation (\ref{s2}), for the Poisson bracket $\{ I, J \}_1$,
we have
$$
{\pa g^{ij}_1 \over \pa u^k} = {\pa P^{ij} \over \pa u^k}
+ {\pa P^{ji} \over \pa u^k} + K_1 \delta^j_k u^i + K_1 \delta^i_k u^j
$$
and, consequently, taking into account relation (\ref{s1}),
we get
$$
g^{ij}_1 (u) = P^{ij} (u) + P^{ji} (u) + K_1 u^i u^j + c^{ij},
\ \ \ c^{ij} = {\rm const}, \ \ \ c^{ij} = c^{ji}.
$$
Thus it is proved that the coefficients of the Poisson bracket
$\{ I, J \}_1$ have the so-called Liouville form
(see about the important Liouville property more in detail below):
$$
g^{ij}_1 (u) = R^{ij} (u) + R^{ji} (u) + K_1 u^i u^j,
$$
$$
b^{ij}_{1, k} (u) = {\pa R^{ij} \over \pa u^k} + K_1 u^i \delta^j_k,
$$
where
$$
R^{ij} (u) = P^{ij} (u) + {1 \over 2} c^{ij}.
$$

Moreover, from relation (\ref{c1}) we get additionally
$$
\eta_{pj} b^{jr}_{1, l} (u) = \eta_{li} b^{ir}_{1, p} (u),
$$
that is,
$$
{\pa (\eta_{ps} R^{sr} (u)) \over \pa u^l} +
K_1 \eta_{ps} u^s \delta^r_l =
{\pa (\eta_{ls} R^{sr} (u)) \over \pa u^p} +
K_1 \eta_{ls} u^s \delta^r_p.
$$
The last formula is equivalent to the relation
$$
{\pa (\eta_{ps} R^{sr} (u) +
K_1 \eta_{ps} u^s u^r ) \over \pa u^l} =
{\pa (\eta_{ls} R^{sr} (u) +
K_1 \eta_{ls} u^s u^r) \over \pa u^p} .
$$
Consequently, by the Poincar\'e lemma,
there locally exist functions
$H^r (u)$ such that
$$
\eta_{ps} R^{sr} (u) + K_1 \eta_{ps} u^s u^r = {\pa H^r \over
\pa u^p}.
$$
Thus, it is proved that
$$
R^{sr} (u) = \eta^{sp} {\pa H^r \over \pa u^p} - K_1 u^s u^r,
$$
and the coefficients of the Poisson bracket
$\{ I, J \}_1$  have the form:
\be
g^{ij}_1 (u) = \eta^{is} {\pa H^j \over \pa u^s}
+\eta^{js} {\pa H^i \over \pa u^s} - K_1 u^i u^j,
\ee
\be
b^{ij}_{1, k} (u) = \eta^{is} {\pa^2 H^j \over \pa u^s \pa u^k}
- K_1 \delta^i_k u^j.
\ee

Since in this case, as is easy to check, all the
relations of compatibility
(\ref{c1}) and (\ref{c2}) are satisfied, theorem
\ref{t1} is proved.

\section{Integrable equations for canonical \\
compatible pairs of brackets}

\bt \label{liu}
An arbitrary nonlocal bracket of the form
(\ref{nonloc3}) (may be degenerate) is a Poisson bracket
if and only if the following equations are satisfied:
\be
{\pa^2 H^i \over \pa u^k \pa u^s} \eta^{sp}
{\pa^2 H^j \over \pa u^p \pa u^l} =
{\pa^2 H^j \over \pa u^k \pa u^s} \eta^{sp}
{\pa^2 H^i \over \pa u^p \pa u^l}, \label{ass1}
\ee
\bea
&&
\left ( \eta^{ir} {\pa H^s \over \pa u^r} +
\eta^{sr} {\pa H^i \over \pa u^r} - K_1 u^i u^s \right ) \eta^{jp}
{\pa^2 H^k \over \pa u^p \pa u^s} =\nn\\
&&
\left ( \eta^{jr} {\pa H^s \over \pa u^r} +
\eta^{sr} {\pa H^j \over \pa u^r} - K_1 u^j u^s
 \right ) \eta^{ip}
{\pa^2 H^k \over \pa u^p \pa u^s}. \label{ass2}
\eea
\et

In the flat case ($K_1 = 0$),
the corresponding theorem was stated by the present
author in \cite{16}, where was also stated
the conjecture on the integrability of the corresponding
system (\ref{ass1}),
(\ref{ass2}) (see also the conditions on the flat pencils
of metrics in \cite{18}).

For every nonlocal bracket of the form
(\ref{nonloc3}) relations (\ref{s1}), (\ref{s2}), (\ref{s5})
are satisfied identically.
Relation (\ref{s4}) takes the form
\be
b^{ij}_{1, s} (u) b^{sr}_{1, k} (u) =
b^{ir}_{1, s} (u) b^{sj}_{1, k} (u)  \label{assoc}
\ee
that gives the equations
\bea
&&
\left ( \eta^{ip} {\pa^2 H^j \over \pa u^p \pa u^s}
- K_1 \delta^i_s u^j \right )
\left ( \eta^{sl} {\pa^2 H^r \over \pa u^l \pa u^k}
- K_1 \delta^s_k u^r \right ) = \nn\\
&&
\left ( \eta^{ip} {\pa^2 H^r \over \pa u^p \pa u^s}
- K_1 \delta^i_s u^r \right )
\left ( \eta^{sl} {\pa^2 H^j \over \pa u^l \pa u^k}
- K_1 \delta^s_k u^j \right )
\eea
equivalent to the equations (\ref{ass1}).

Relation (\ref{s3}) gives the equations
\bea
&&
\left ( \eta^{ip}  {\pa H^s \over \pa u^p} +
\eta^{sp} {\pa H^i \over \pa u^p} - K_1 u^i u^s \right )
\left ( \eta^{jl} {\pa^2 H^r \over \pa u^l \pa u^s} -
K_1 \delta^j_s u^r \right ) =\nn\\
&&
\left ( \eta^{jp}  {\pa H^s \over \pa u^p} +
\eta^{sp} {\pa H^j \over \pa u^p} - K_1 u^j u^s \right )
\left ( \eta^{il} {\pa^2 H^r \over \pa u^l \pa u^s} -
K_1 \delta^i_s u^r \right )
\eea
equivalent to the equations (\ref{ass2}).

\bc
Any linear function $H^i = c^i_k u^k + c^i,$
$c^i_k = {\rm const},$ $c^i = {\rm const},$
is a trivial solution of the
nonlinear system of equations
(\ref{ass1}), (\ref{ass2}). Thus the bracket
\bea
&&
\ \{ I,J \}_1 = \int {\delta I \over \delta u^i(x)}
\left ( \left [ \eta^{is} c^j_s + \eta^{js} c^i_s
- K_1 u^i u^j \right ] {d \over dx} - \right. \nn\\
&&
\left. K_1 \delta^i_k u^j
\, u^k_x
+ K_1 u^i_x \left ( {d \over dx} \right )^{-1} u^j_x \right )
{\delta J \over \delta u^j(x)} dx
\label{nonloc3aa}
\eea
is a Poisson bracket for any constants $c^i_k$.
In particular, for any symmetric constant matrix
$(\mu^{ij}),$ $\mu^{ij} = \mu^{ji},$
$\mu^{ij} = {\rm const},$ the contravariant metric
$\mu^{ij} - K u^i u^j$ is always a metric of constant
Riemannian curvature
$K$ if it is nondegenerate (it is obvious
that under the condition
of nondegeneracy for the matrix $(\mu^{ij})$ this metric
is nondegenerate at least for small $u^i$),
and in addition the Levi-Civita connection
generated by this metric is defined by the relation
$\Gamma^{ij}_k (u) =
(\mu^{is} - K u^i u^s) \Gamma^j_{sk} (u)
=  K \delta^i_k u^j.$
\ec

Consider the contravariant metric of the form
\be
g^{ij} (u) =
a^i \delta^{ij} - K u^i u^j, \label{can}
\ee
where $a^i,$ $1 \leq i \leq N,$ are
arbitrary nonzero constants, so that the metric
$g^{ij} (u)$ is nondegenerate.
It is easy to prove that
\be
\det (g^{ij} (u)) = a^1 \cdots a^N \left (1 - K \sum_{s=1}^N
{(u^s)^2 \over a^s} \right ).
\ee
The covariant metric $g_{ij} (u),$ which is inverse to the
contravariant metric $g^{ij} (u),$  $g^{is} (u) g_{sj} (u) =
\delta^i_j,$ has the form
\be
g_{ij} (u) =
{1 \over a^i} \delta_{ij} + {K u^i u^j
\over a^i a^j \left ( 1 - K \sum_{s=1}^N {(u^s)^2 \over a^s} \right )}
\ee
and is a metric of constant Riemannian
curvature $K$ for all values of the nonzero
constants $a^i$. The considered local
coordinates are geodesic at the point with the
coordinates $(0,...,0)$ ($u^i = 0$), but
they are not the normal Birkhoff coordinates and
even not Riemannian coordinates.

Note that the metric $g^{ij} (u)$ is nondegenerate
if and only if not more than one of
$N + 1$ constants: $a^i,$ $1 \leq i \leq N,$ and $K$, equals to
zero. All the formulae are easily adapted for the case
when one of these constants equals to zero.
If $a^m = 0$, then
\be
\det (g^{ij} (u)) = - K \left ( \prod_{s \neq m} a^s \right ) (u^m)^2.
\ee
In this case, the components of the covariant metric
$g_{ij} (u)$  of constant Riemannian curvature $K$
have the form
\be
g_{mm} (u) = - {1 \over K (u^m)^2}  \left (
1 - K \sum_{s \neq m} {(u^s)^2 \over
a^s} \right ),
\ee
\be
g_{im} (u) = g_{mi} (u) = - {1 \over a^i} {u^i \over u^m},\ \ \
i \neq m,
\ee
\be
g_{ii} (u) = { 1 \over a^i}, \ \ \ i \neq m,
\ee
\be
g_{ij} (u) = 0, \ \ \ i \neq j, \  i \neq m, \
j \neq m.
\ee
All these models of the spaces
of constant curvature play an important role
in the Hamiltonian theory of systems of hydrodynamic type.
The Mokhov--Ferapontov brackets generated
by the metrics of constant Riemannian curvature
(\ref{can}) for
$a^i = \varepsilon^i = \pm 1$ are called in \cite{19}
{\it canonical}. The canonical brackets arose naturally
also in applications in
\cite{20}.

\bt [\cite{2}, \cite{3}] \label{int}
The system of nonlinear equations (\ref{ass1}), (\ref{ass2})
is integrated by the method of inverse scattering
problem.
\et

Note that in \cite{2}, \cite{3}
the system of nonlinear equations describing the
compatible nonlocal Poisson brackets of the form
(\ref{nonl}) in other local coordinates, more
convenient for the integration (in these coordinates, the
metrics of both compatible brackets are diagonal),
was derived and integrated.

\section{Liouville and special Liouville coordinates}

Local coordinates $u = (u^1,..., u^N)$ are called
{\it Liouville} for an arbitrary Poisson bracket
$\{ I, J \}$ if the functions (the fields) $u^i (x)$
are densities of integrals in involution with
respect to this bracket, that is,
\be
\ \{ U^i, U^j \} = 0, \ \ \ 1 \leq i, j \leq N,
\ee
where $U^i = \int u^i (x) dx,$ $1 \leq i \leq N.$
In this case the Poisson bracket is also
called {\it Liouville} in these coordinates.
Liouville coordinates naturally arise and
play an essential role in the Dubrovin--Novikov
procedure of averaging of Hamiltonian equations \cite{8}.
The physical coordinates derived by averaging
of the densities of the participating in the Dubrovin--Novikov
procedure $N$
involutive local integrals of an initial
Hamiltonian system are always Liouville for the
corresponding averaged bracket.
This property was a motivation for
the definition of Liouville coordinates
for local Poisson brackets of hydrodynamic type in \cite{8}.
For general nonlocal Poisson brackets of hydrodynamic type
(the Ferapontov brackets \cite{21}, \cite{22}),
Liouville coordinates were introduced in \cite{19}.

A nonlocal Poisson bracket of hydrodynamic type
(\ref{nonl}) is Liouville in the local coordinates
$u = (u^1,...,u^N)$ if and only if there exists
a matrix function $\Phi^{ij} (u)$ such that
\be
b^{ij}_k (u) = {\pa \Phi^{ij} \over \pa u^k}
- K \delta^i_k u^j.
\ee
In this case, by virtue of relations
(\ref{s1}), (\ref{s2}), the metric $g^{ij} (u)$
must have the form
\be
g^{ij} (u) = \Phi^{ij} (u) + \Phi^{ji} (u) -
K u^i u^j
\ee
(here the function $\Phi^{ij} (u)$ can be corrected
by a constant matrix function
$c^{ij} = {\rm const}$).
The matrix function $\Phi^{ij} (u)$ is called {\it a Liouville
function}.

Thus a nonlocal Poisson bracket (\ref{nonl}) is
Liouville if it has the form
\bea
&&
\ \{ I,J \} = \int {\delta I \over \delta u^i(x)}
\left ( \left ( \Phi^{ij} (u) + \Phi^{ji} (u) - K u^i u^j \right )
{d \over dx} + \right.  \nn\\
&&
\ \left. \left ( {\pa \Phi^{ij} \over
\pa u^k} - K \delta^i_k u^j  \right ) \, u^k_x
+ K u^i_x \left ( {d \over dx} \right )^{-1} u^j_x \right )
{\delta J \over \delta u^j(x)} dx.
\label{nonlliu}
\eea

From theorem \ref{t1}, it follows

\bt
Flat coordinates of an arbitrary nondegenerate
local Poisson bracket of hydrodynamic type $\{ I, J \}_2$
are always Liouville for any nonlocal Poisson bracket
$\{ I, J \}_1$ of the form (\ref{nonl}) compatible with
$\{ I, J \}_2$.
Moreover, in addition the corresponding Liouville
function $\Phi^{ij} (u)$
always has the special form
\be
\Phi^{ij} (u) = \eta^{is} {\pa H^j \over \pa u^s}.
\ee
\et

Local coordinates $u = (u^1,...,u^N)$ are called
{\it special Liouville coordinates} \cite{16},
\cite{17}
for an arbitrary Poisson bracket
$\{ I, J \}$ if there exists a nonzero constant
symmetric matrix $(\eta_{ij})$ such that
the functions (the fields) $u^i (x),$ $1 \leq i \leq N,$
and $\eta_{ij} u^i (x) u^j (x)$ are densities of
integrals in involution with respect to
this bracket, that is,
\be
\ \{ U^i, U^j \} = 0, \ \ \ 1 \leq i, j \leq N + 1,
\ee
where $U^i = \int u^i (x) dx,$ $1 \leq i \leq N,$
$U^{N+1} = \int \eta_{ij} u^i (x) u^j (x) dx$.
In this case the Poisson bracket is also called
{\it special Liouville} in these coordinates.
The special Liouville coordinates were introduced in
\cite{16}, \cite{17}. The most important case
is the case of nondegenerate matrix $\eta_{ij}$.

\bt
An arbitrary nonlocal Poisson bracket
of the form (\ref{nonl}) is special Liouville
in the local coordinates
$u = (u^1,..., u^N)$ if and only if
it is Liouville with a special Liouville function
$\Phi^{ij} (u)$ such that
\be
\eta_{ks} \Phi^{sj} (u) = {\pa H^j \over \pa u^k}.
\ee
\et

In this case, for a nondegenerate matrix $(\eta_{ij})$,
we get exactly our bracket (\ref{nonloc3})
from the canonical compatible pair.

Thus our problem on compatible nonlocal
Poisson brackets of hydrodynamic type is
equivalent to the
problem of classification of the
special Liouville coordinates for
nonlocal Poisson brackets of hydrodynamic type.

\bt
An arbitrary nonlocal Poisson bracket of hydrodynamic type
of the form (\ref{nonl}) is compatible
with the constant Poisson bracket (\ref{lokconst})
if and only if the functions $u^i (x),$ $1 \leq i \leq N,$
and $\eta_{ij} u^i (x) u^j (x),$ $\eta^{is} \eta_{sj} =
\delta^i_j,$ are densities of integrals in
involution with respect to the Poisson bracket (\ref{nonl}).
\et

Note that $u^i (x),$ $1 \leq i \leq N,$
are the densities of the annihilators of the bracket
(\ref{lokconst}), and ${1 \over 2} \eta_{ij}
u^i (x) u^j (x)$ is the density of the momentum
of the bracket (\ref{lokconst}).

\bt
An arbitrary nonlocal Poisson bracket of
hydrodynamic type of the form (\ref{nonl}) is compatible
with an arbitrary nondegenerate local Poisson bracket
of hydrodynamic type
(\ref{lok}) if and only if $N$ annihilators and the momentum
of the bracket (\ref{lok}) are integrals in involution
with respect to the Poisson bracket (\ref{nonl}).
\et

\section{Integrable bi-Hamiltonian systems\\
of hydrodynamic type}

Consider an arbitrary pair of compatible
nonlocal Hamiltonian
operators of hydrodynamic type
$P^{ij}_1$ and $P^{ij}_2$
generated by metrics of constant Riemannian curvature.
As is shown above in lemma \ref{locnon},
one of these operators, let us assume $P^{ij}_2$,
can be considered as local without
loss of generality. If the local Hamiltonian operator
$P^{ij}_2$ is nondegenerate, then it follows
from theorem \ref{t1} that, by local change
of coordinates, the pair of compatible Hamiltonian
operators
$P^{ij}_1$ and $P^{ij}_2$ can be reduced to the following
canonical form:

\be
P^{ij}_2 [v(x)] = \eta^{ij} {d \over dx}, \label{op1}
\ee
\be
P^{ij}_1 [v(x)] =
\left ( \eta^{is} {\partial h^j \over \partial v^s}
+ \eta^{js} {\partial h^i \over \partial v^s} - K v^i v^j \right )
{d \over dx} + \left ( \eta^{is}
{\partial^2 h^j \over \partial v^s \partial v^k} -
K \delta^i_k  v^j \right )
v^k_x + K v^i_x \left ( {d \over dx} \right )^{-1} v^j_x,   \label{op2}
\ee
where ($\eta^{ij}$) is an arbitrary nondegenerate
constant symmetric matrix:
$\det (\eta^{ij}) \neq 0,$ $ \eta^{ij} = {\rm const},$
$\eta^{ij} = \eta^{ji};$ $K$ is an arbitrary constant;
$h^i (v),$  $1 \leq i \leq N,$
are smooth functions defined in a certain domain of
local coordinates and such that the operator
(\ref{op2}) is Hamiltonian, that is, the functions
$h^i (v)$ satisfy to the integrable equations
(\ref{ass1}), (\ref{ass2})
(see theorems \ref{liu} and \ref{int} above).

\br
{\rm It is obvious that here we can always
consider that $\eta^{ij} = \varepsilon^i \delta^{ij},$
$\varepsilon^i = 1$ for $i \leq p,$
$\varepsilon^i = - 1$ for $i > p,$
where $p$ is the positive index of inertia of the metric,
$0 \leq p \leq N$, and, in addition,
it is necessary to classify the Hamiltonian operators
(\ref{op2}) with respect to the action of
the group of motions for the corresponding
$N$-dimensional pseudo-Euclidean space ${\bf R}^N_p$,
but for our purposes it is sufficient
(and more convenient) to use the indicated above
representation for canonical compatible pair
(``conventionally canonical''
representation).}
\er

Consider the recursion operator generated by
canonical compatible Hamiltonian operators (\ref{op1}),
(\ref{op2}):

\bea
&&
R^i_l [v(x)] =
\left [ P_1 [v (x)] \left ( P_2 [v(x)] \right )^{-1} \right ]^i_l
=
\left ( \biggl ( \eta^{is} {\partial h^j \over \partial v^s}
+ \eta^{js} {\partial h^i \over \partial v^s}
- K v^i v^j \biggr )
{d \over dx} + \right. \nn\\
&&
\left. \left (\eta^{is}
{\partial^2 h^j \over \partial v^s \partial v^k}
- K \delta^i_k v^j \right )
v^k_x + K v^i_x \left (
{d \over dx} \right )^{-1} v^j_x
\right ) \eta_{jl} \left ( {d \over dx} \right )^{-1}
\label{recur}
\eea
(what about recursion operators generated
by pairs of compatible Hamiltonian operators, see
\cite{23}--\cite{27}, \cite{13}).

Let us apply the derived recursion operator
(\ref{recur}) to the system of translations with respect to
$x$, that is, the system of hydrodynamic type
\be
v^i_t = v^i_x,  \label{transl}
\ee
which is, obviously, Hamiltonian with
the Hamiltonian operator (\ref{op1}):
\be
v^i_t =v^i_x \equiv P^{ij}_2 {\delta H \over \delta v^j (x)},
\ \ \ H = {1 \over 2} \int \eta_{jl} v^j (x) v^l (x) dx.
\ee

Any system from the hierarchy
\be
v^i_{t_n} = \left ( R^n \right )^i_j v^j_x,
\ \ \ \ n \in {\bf Z}, \label{ierarkh}
\ee
is a multi-Hamiltonian integrable system.

In particular, any system of the form
\be
v^i_{t_1} = R^i_j v^j_x,
\ee
that is, the system of hydrodynamic type
\bea
&&
v^i_{t_1} =
\left ( \biggl ( \eta^{is} {\partial h^j \over \partial v^s}
+ \eta^{js} {\partial h^i \over \partial v^s} -
K v^i v^j \biggr )
{d \over dx} + \left (\eta^{is}
{\partial^2 h^j \over \partial v^s \partial v^k} -
K \delta^i_k v^j \right )
v^k_x + \right. \nn\\
&&
\left. K v^i_x \left (
{d \over dx} \right )^{-1} v^j_x \right ) \eta_{jl} v^l
\equiv  \left (  \eta^{is} {\partial h^j \over \partial v^s} \eta_{jk}
+  {\partial h^i \over \partial v^k}
 + \eta^{is} \eta_{jl}
{\partial^2 h^j \over \partial v^s \partial v^k} v^l
- K \eta_{sk} v^i v^s - \right. \nn\\
&&
\left. {K \over 2} \delta^i_k
\eta_{sl} v^s v^l \right )
v^k_x
\equiv
\left ( h^i (v)+ \eta^{is} {\partial h^j \over
\partial v^s} \eta_{jl} v^l - {K \over 2}
\eta_{sk} v^i v^s v^k \right )_x,  \label{canon}
\eea
where $h^i(v),$ $1 \leq i \leq N,$ is an arbitrary
solution of the integrable system (\ref{ass1}), (\ref{ass2}),
is integrable.

This system of hydrodynamic type is bi-Hamiltonian
with the pair of canonical compatible Hamiltonian operators
(\ref{op1}),
(\ref{op2}):
\bea
&&
v^i_{t_1} =
\left ( \biggl ( \eta^{is} {\partial h^j \over \partial v^s}
+ \eta^{js} {\partial h^i \over \partial v^s} -
K v^i v^j \biggr )
{d \over dx} + \left (\eta^{is}
{\partial^2 h^j \over \partial v^s \partial v^k} -
K \delta^i_k v^j \right )
v^k_x + \right.\nn\\
&&
\left. K v^i_x \left (
{d \over dx} \right )^{-1} v^j_x \right )
{\delta H_1 \over \delta v^j (x)}, \ \ \
H_1 = {1 \over 2} \int \eta_{jl} v^j (x) v^l (x) dx, \label{eq1}
\eea
\be
v^i_{t_1} = \eta^{ij} {d \over dx}
{\delta H_2 \over \delta v^j (x)}, \ \ \
H_2 =  \int \left (\eta_{jk} h^k (v (x)) v^j (x)
- {K \over 8} \eta_{jk} \eta_{sl} v^j v^k v^s v^l \right )
dx. \label{eq2}
\ee

The next system in the hierarchy (\ref{ierarkh}) (for $n=2$)
is the integrable system of hydrodynamic type
\bea
&&
v^i_{t_2} =
\left ( \biggl ( \eta^{is} {\partial h^j \over \partial v^s}
+ \eta^{js} {\partial h^i \over \partial v^s} -
K v^i v^j \biggr )
{d \over dx} + \left (\eta^{is}
{\partial^2 h^j \over \partial v^s \partial v^k} -
K \delta^i_k v^j \right )
v^k_x + \right. \nn\\
&&
\left. K v^i_x \left (
{d \over dx} \right )^{-1} v^j_x \right )
\eta_{jl}
\left ( h^l (v)+ \eta^{lp} {\partial h^r \over
\partial v^p} \eta_{rq} v^q
- {K \over 2}
\eta_{pr} v^l v^p v^r
\right )
\equiv \nn\\
&&
\left ( \biggl ( \eta^{is} {\partial h^j \over \partial v^s}
+ \eta^{js} {\partial h^i \over \partial v^s}
- K v^i v^j \biggr )
\left ( \eta_{jl} {\partial h^l \over \partial v^k} +
\eta_{rk} {\partial h^r \over \partial v^j} +
\eta_{rq} v^q {\partial^2 h^r \over \partial v^j
\partial v^k} - \right. \right.
\nn\\
&&
 \left. K \eta_{jl} \eta_{pk} v^l v^p -
{K \over 2} \eta_{jk} \eta_{pl} v^l v^p \right ) +
\left ( \eta^{is}
{\partial^2 h^j \over \partial v^s \partial v^k} -
K \delta^i_k v^j  \right )
\left ( \eta_{jl} h^l (v) + \right.\nn\\
&&
\left. \left.
\eta_{rq} v^q {\partial h^r \over \partial v^j} -
{K \over 2} \eta_{jl} \eta_{pr} v^l v^p v^r \right )
+ K \delta^i_k  \left ( \eta_{jl} h^l (v) v^j -
 {K \over 8} \eta_{jl} \eta_{pr}
v^j v^l v^p v^r \right )
\right )
v^k_x.
\eea

Even the trivial, linear with respect to the fields $v^i (x)$,
solutions of the system (\ref{ass1}), (\ref{ass2})
generate nontrivial integrable systems of
hydrodynamic type.
All the results presented in this work are generalized
directly to the important case of general nonlocal Poisson
brackets of hydrodynamic type (the Ferapontov brackets
 \cite{21}, \cite{22}),
although the corresponding formulae become
considerably more cumbersome and less
effective. These results will be published
in another paper.

\begin{flushleft}
Centre for Nonlinear Studies,\\
L.D.Landau Institute for Theoretical Physics, \\
Russian Academy of Sciences\\
e-mail: mokhov@mi.ras.ru; mokhov@landau.ac.ru\\
\end{flushleft}

\end{document}